\documentclass[12pt]{article}
\usepackage{hyperref}
\usepackage{graphicx,amsmath,amsfonts,latexsym,amssymb,amsthm,amscd}
\usepackage{latexsym,amscd}

  \newcommand{\field}[1]{\mathbb{#1}}
  \newcommand{\rbb}{\field{R}}
  \newcommand{\cbb}{\field{C}}

  \newcommand{\I}{\,\mathrm{i}\,}
  \newcommand{\PsDO}{\Psi\mathrm{DO}}

  \newcommand{\tr}{\mathrm{Tr}}

  \newcommand{\beq}{\begin{equation}}
  \newcommand{\eeq}{\end{equation}}

  \newtheorem{theorem}{Theorem}[section]

  \newtheorem{exa}[theorem]{Example}

\begin{document}

\title{Geometry of the high energy limit of differential operators on vector bundles}

\author{Alexander Strohmaier}

\maketitle

\begin{abstract}
 At high energies relativistic quantum systems describing scalar particles behave classically. This observation
 plays an important role in the investigation of eigenfunctions of the Laplace operator on manifolds
 for large energies and allows to establish relations to the dynamics
 of the corresponding classical system. Relativistic quantum systems describing particles
 with spin such as the Dirac equation do not behave classically at high energies. 
 Nonetheless, the dynamical properties of the classical frame flow determine the behavior of 
 eigensections of the corresponding operator for large energies.
 We review what a high energy limit is and how it can be described for geometric operators.
\end{abstract}


\section{Introduction}

A quantum physical system is usually described by an algebra $\mathcal{A}$ 
of operators on a Hilbert space $\mathcal{H}$ and the time evolution, which is
a one-parameter group $U(t)$ of unitary operators on $\mathcal{H}$.
An important example of such a system is the one describing the motion of a quantum
particle on a compact Riemannian manifold $M$. In this case the Hilbert space is
$L^2(M,\mu_g)$, where $\mu_g$ is the Riemannian measure.
The time evolution is described by
the Schr\"odinger equation, 
which means that the unitary one-parameter group $U(t)$ is given by
$$
 U(t) = e^{-\I t H},
$$
where $H$ is the Schr\"odinger operator. For non-relativistic quantum systems
$H$ would typically be $\Delta + V$, where $\Delta$ is the metric Laplace operator on $M$
and $V \in C^\infty(M,\mathbb{R})$ is a potential. As we are interested in the
high energy limit we will consider the Klein-Gordon time evolution which
describes relativistic particles. So we take
$$
 H= (\Delta + m^2 +V)^{1/2}
$$ 
where $m$ is a positive real number, the mass of the particle. We will assume
here that $m^2+V(x) \geq 0$ so that $H$ is a positive first order pseudodifferential
operator with principal symbol $\sigma_H(\xi)=\Vert \xi \Vert$.
The above group is defined by spectral calculus as
$$
 \Delta + m^2+ V : H^2(M) \to L^2(M)
$$
is a self-adjoint operator with the  Sobolev space $H^2(M)$ as its domain.
The algebra of observables would be a unital $*$-subalgebra of $\mathcal{B}(\mathcal{H})$, the algebra
of all bounded operators on $\mathcal{H}$.
A state on $\mathcal{A}$ is a linear functional $\omega: \mathcal{A} \to \mathbb{C}$
such that
\begin{itemize}
 \item[(i)] $\omega$ is complex linear,
 \item[(ii)] $\omega$ is positive, i.e. $\omega(A^* A) \geq 0$ for all $A \in \mathcal{A}$,
 \item[(iii)] $\omega(1)=1$.
\end{itemize}
The physical interpretation is that a state is a state of the system and
assign to each observable its expectation value.
An example is
$$
 \omega(A) = \langle \psi, A \psi \rangle,
$$
where $\psi$ is a vector of unit length in $\mathcal{H}$.
More generally, any trace class operator $\rho$ with $\rho \geq 0$ and
$\tr(\rho)=1$ defines a state by
$$
 \omega_\rho(A) = \tr(A \rho).
$$
These states are often referred to as normal.
 If $A$ is a self-adjoint element in $\mathcal{A}$
then $A$ has a spectral decomposition
$$
 A = \int_\mathbb{R} \lambda dE_\lambda.
$$
If the system is in the normal state $\omega_\rho$
then the probability that a measurement of the
observable $A$ yields a value in the Borel measurable set 
$\mathcal{O} \subset \mathbb{R}$ is given by
$$
 \omega_\rho(E_{\mathcal{O}}) = \mathrm{tr}(\rho E_{\mathcal{O}}),
$$
where $E_{\mathcal{O}}$ is the spectral projection onto $\mathcal{O}$
$$
 E_{\mathcal{O}} =  \int_\mathbb{\mathcal{O}} dE_\lambda.
$$
This is the probabilistic interpretation of quantum mechanics.
For example if $\mathcal{U}$ is a subset of $M$ then the characteristic function
$\chi_\mathcal{U}$ corresponds to the binary experiment that measures whether the particle is in the region
$\mathcal{U}$. The probability of finding the particle in $\mathcal{U}$ is therefore given
by $\omega_\rho(\chi_\mathcal{U})$
or if the state is of the form $\omega_\psi$ with $\psi \in L^2(M)$ then of course we get for the probability
of finding the particle in $\mathcal{U}$
$$
 \int_\mathcal{U} \vert \psi(x) \vert^2 dx.
$$

Whereas in the literature $\mathcal{B}(\mathcal{H})$ itself is often chosen as the algebra of observables, this choice is sometimes not very convenient for practical purposes. 
On the physical side it is impossible to build a detector that measures $\chi_\mathcal{O}$. The reason is that such a measurement would involve a detector that near the boundary of $\mathcal{O}$
had an arbitrary high resolution. If we wanted to be more realistic we would restrict ourself to algebras that contain functions that are only smooth or continuous.
On the mathematical side it is much easier to specify a state on a smaller algebra rather than on the
full algebra of bounded operators. Knowledge of the state on the smaller algebra is often sufficient
to extend it uniquely to a larger subalgebra. So which algebra to choose for the particle on the manifold?

The state of a classical particle is completely determined by its momentum and its position. We expect the same to be true for
quantum particles. In order to measure the position we take the algebra
$C^\infty(M)$. Measurement of the momentum involves unbounded operators of the form $\I X$, where
$X$ is a vector field. In other words we would need an algebra of operators that contains enough
bounded functions of $X$ so that we can approximate the spectral projections of $\I X$
by elements in our algebra. If we choose a classical symbol $p$ of order $0$ on $\mathbb{R}$
then $p(X)$ is a classical pseudodifferential operator of order $0$. So if we choose the algebra of pseudodifferential operators $\PsDO^0_{cl}(M)$ on $M$ this algebra contains enough observables to measure the location and the momentum of our particle up to some arbitrary small error.
And, indeed, the restriction of a normal state $\omega_\rho$ to the algebra $\mathcal{A}$
determines this state completely. To see this note that normal states are continuous in
the weak-$*$-topology on $\mathcal{B}(\mathcal{H})$ and $\PsDO^0_{cl}(M)$
is weak-$*$-dense in $\mathcal{B}(\mathcal{H})$. 
Since states are automatically norm continuous any state can by continuity be uniquely extended to a state on the norm closure it is reasonable to use as the algebra of observables the norm closure
of the algebra of pseudodifferential operators, that is
$$
 \mathcal{A} = \overline{\PsDO^0_{cl}(M)}.
$$

The high energy limit is now the limit of the quantum mechanical system for particles of high energies.
Measuring the energy of a particle corresponds to the unbounded operator
$$
 \Delta + V
$$
which generates the unitary one-parameter group $U(t)$. So we would think of a normal state $\omega_\rho$ as a state with energy larger than $\lambda$ if
$$
 \omega_\rho (E_{[-\infty,\lambda)})=0
$$
where $E$ is the spectral projection of $ \Delta + V$. 

Note that the set of states is weak-$*$-compact and the algebra $\mathcal{A}$ is separable.
Therefore, the set of states is sequentially compact in the weak-$*$-topology. That is every
sequence of normal states has weak-$*$ limit points. Suppose
$\omega_{\rho_n}$
is a sequence of normal states that converges to a not necessarily normal state $\omega_\infty$
in the weak-$*$-topology. Then we think of $\omega_\infty$ as a high energy limit
if
$$
 \omega_{\rho_n}(E_{[-\infty,\lambda)}) \to  0
$$
for all $\lambda>0$.
Using spectral calculus one finds that this is equivalent to
$$
 \omega_\infty((\vert \Delta+V \vert+1)^{-1}) =0.
$$
 Now suppose that $K$ is a pseudodifferential operator of order $-1$.
Then, $K \sqrt{\vert \Delta+V \vert+1}$ is a pseudodifferential operator
of order $0$ and is therefore bounded. By the Cauchy-Schwarz inequality
$$
 \vert \omega_{\infty}(K) \vert^2 \leq \omega_{\infty}\left(\left\vert K \sqrt{\vert \Delta+V \vert+1} \right\vert^2 \right) \cdot
 \omega_{\infty}\left((\sqrt{\vert \Delta+V \vert+1})^{-2} \right) 
$$
and therefore
$$
  \omega_{\infty}(K)  = 0
$$
for any high energy limit $\omega_\infty$. Therefore, $\omega_\infty$ vanishes on the algebra
of pseudodifferential operators of order $-1$ and since it is continuous it also vanishes on
the ideal of compact operators $\mathcal{K}$. Consequently, at high energies the states become states on the quotient algebra
$$
 \mathcal{A} / \mathcal{K}.
$$
Now it is well known (e.g. \cite{MR0173174}, Th. 11.1) that the principal symbol map
$$
 \sigma:  \PsDO^0_{cl}(M) \to C^\infty(S^* M)
$$
extends continuously to a map
$$
 \hat \sigma: \mathcal{A} \to C(S^* M)
$$
and $\ker \hat \sigma = \mathcal{K}$.
This means that
$$
 \mathcal{A} / \mathcal{K}
$$
is naturally isomorphic to the commutative algebra $C(S^* M)$.
States on $C(S^* M)$ are by the Riesz representation theorem in one to one
correspondence to regular Borel probability measures on $S^* M$. So for every high energy limit
$\omega_\infty$ there exists a unique probability measure $\mu$ on $S^* M$ such that
$$
 \omega_\infty(A) = \int_{S^* M} \hat \sigma_A(\xi) d\mu(\xi).
$$
The fact that the high energy limit states are actually states on an abelian algebra
can be interpreted as the passage from quantum to classical mechanics. The system behaves
classically for very large energies. That the quantum mechanical time evolution becomes the
classical motion along geodesics is now a consequence of Egorov's theorem.
Namely, if $A \in \mathcal{A}$ then also
$$
 A(t)=U(-t) A U(t) \in \mathcal{A}
$$
and
$$
 \hat \sigma_{A(t)} = G_t^* ( \hat\sigma_{A}),
$$ 
where $G_t$ is the geodesic flow on $S^*M$ and $G_t^*$ is its pull-back
acting on functions. This means that the group of $*$-automorphisms
$$
 \alpha_t(A)=U(-t) A U(t)
$$
which describes the quantum mechanical time evolution on the level of observables
(the so-called Heisenberg picture) factors to $\mathcal{A}/\mathcal{K} \cong C(S^* M)$
and becomes there the geodesic flow. This is a very concise way of saying that in the limit
of high energy the quantum system becomes classical and the time-evolution becomes
the motion along geodesics with constant speed.

Note that in the high energy limit the potential does not play a role any more. From the physical
point of view this is expected as particles with high energy do not  ''feel'' a potential and move at the
speed of light along lightlike geodesics.

Interestingly some high energy limit can be computed explicitly.
Let $\phi_j$ be a complete orthonormal sequence of eigenfunctions of $\Delta$ such that
\begin{gather*}
 \Delta \phi_j = \lambda_j^2 \phi_j,\\
 0=\lambda_1 \leq \lambda_2 \leq \lambda_3 \leq \ldots. 
\end{gather*}
Then, the sequence of normal states $\omega_N$ defined by
$$
 \omega_N(A) = \frac{1}{N} \sum_{j=1}^N \langle \phi_j , A \phi_j \rangle
$$
has a high energy limit. It follows from the classical Tauberian theorem of Karamata
that the limit of this sequence is given by
$$
 \omega_\infty(A)=\lim_{t \to 0^+} \frac{\tr (A e^{-t (\Delta+V)})}{\tr(e^{-t (\Delta+V)})}.
$$
The (microlocal) heat kernel expansion then shows that
$$
  \lim_{t \to 0^+} \frac{\tr (A e^{-t (\Delta+V)})}{\tr(e^{-t (\Delta+V)})} = \int_{S^* } \hat \sigma_A(\xi) d\mu_L(\xi),
$$
where $\mu_L$ is the normalized Liouville measure on $S^* M$.

The state $\omega_t$ defined by
$$
 \omega_t(A)=\frac{\tr (A e^{-t (\Delta+V)})}{\tr(e^{-t (\Delta+V)})}
$$
is the
KMS-state with temperature $t^{-1}$ describing a quantum system at temperature
$t^{-1}$ in thermal equilibrium. In the limit as the temperature goes to
infinity the state converges to the Liouville measure on the unit-cotangent bundle.

The sequence of eigenstates $\langle \phi_j, \cdot \phi_j \rangle$ is a sequence of invariant
states and any weak-$*$-limit point is a therefore an invariant high energy limit. The above says that on average these
states converge to the Liouville measure. If the Liouville measure is ergodic with respect to the
geodesic flow this means that the tracial state $\omega_\infty$ defined above is ergodic. This means there is no non-trivial decomposition of $\omega_\infty$ into a convex combination of invariant states.
From this one can conclude that any subsequence of $\langle \phi_j, \cdot \phi_j \rangle$
that does not have the state $\omega_\infty$ as a weak-$*$-limit has to have counting density
zero, as otherwise it would give rise to a decomposition of $\omega_\infty$ into invariant states.
A more careful argument along these lines (see
\cite{MR0402834, MR1239173, MR0818831, MR0916129})
shows that in fact there is a subsequence of counting density one
of eigenfunctions $\phi_{j(k)}$ such that
$$
 \lim_{k \to \infty} \langle \phi_{j(k)}, A \phi_{j(k)} \rangle = \omega_{\infty}(A),
$$
for all $A \in \mathcal{A}$. This is usually referred to as Quantum ergodicity.

\section{The Dirac equation and Laplace type operators}

Relativistic quantum systems that describe particles with spin
like electrons, neutrinos are not described by the Klein-Gordon equation. 
As the particles have an internal degree of freedom, the spin, 
the Hilbert space will consist of vector valued functions and the observable algebra
needs to include operators that detect these internal degrees of freedom.
This is appropriately described by the following construction.
Suppose that $E \to M$ is a complex hermitian vector bundle over $M$ we take
as the Hilbert space the space of square integrable sections of $E$
$$
 \mathcal{H}=L^2(M;E)
$$
and as an algebra of observables we take the norm closure of the space of zero order classical
pseudodifferential operators acting on sections of this vector bundle
$$
 \mathcal{A}=\overline{\PsDO^0_{cl}(M;E)}.
$$
Now a second order differential operator $P$ is said to be of Laplace-type
if in local coordinates it has the form
$$
 P = - \sum_{ij} g^{ij} \frac{\partial^2}{\partial x^i \partial x^j} + B,
$$
where $B$ has order $1$, or in other words if and only if
$$
 \sigma_P(\xi) = g(\xi,\xi).
$$
Similarly, a first order differential operator $D$ is said to be of Dirac type
if and only if
$$
 \sigma_D(\xi)^2 = g(\xi,\xi).
$$
Of course a first order operator is of Dirac type if and only if its square is
a Laplace type operator.

As we saw in the previous section the Klein-Gordon operator $\Delta + V + m^2$
is a Laplace type operator acting on the trivial vector bundle. So the time
evolution in this case is described by the square root of a Laplace type operator.
To describe electrons one typically chooses a spin
structure on $M$ and then the complex vector bundle is the associated
spinor bundle $S$. The algebra of observables in this case is the algebra of zero order classical pseudodifferential operators
acting on sections of the spinor bundle.
The time-evolution is described by the Dirac
operator $D$ acting on sections of the spinor bundle. This operator
will however not be positive any more, which is a typical feature of relativistic
quantum theory. The system describes electrons and positrons at the same time. Pure electron
states are states that are supported in the positive
spectral subspace of $D$. Since on this subspace the operators
$D$ and $\vert D \vert$ coincide the time evolution of such
states may as well be described by the operator $H=\vert D \vert$.
Whereas from the viewpoint of a one-particle theory it might seem strange to take the operator
$\vert D \vert$ as the generator of the time evolution this is perfectly justified in a fully quantized theory.
The generator of the time-evolution of the fully quantized free electron-positron field restricted
to the one-particle subspace is given by $\vert D \vert$ rather than $D$. The apparent violation
of Einstein causality by the infinite propagation speed of the operator $\mathrm{exp}(-\I t \vert D \vert)$
is resolved in the fully quantized theory and is not causing a problem there (see e.g. \cite{MR1219537} for details).

Note that if one chooses the group generated by $D$ instead of $\vert D \vert$
the time evolution does not leave the space of operators invariant. Instead of passing to
the $\vert D \vert$ one can also restrict the algebra of observables. This approach in favored
in \cite{MR722924} and also used in \cite{MR1644088, MR1694732, MR1838911,MR2073612, BoG04.2} in order to investigate the semi-classical limit
of the Dirac operator.

Spin $1$ particles like photons and mesons are described by Maxwell's
equation or the Proca equation. For example the quantum system describing
photons is given as follows. The vector bundle is the complexified co-tangent bundle
$\Lambda^1 M = T^* M$. The Hilbert space is the closure of the space
of co-closed $1$-forms in $L^2(M;E)$ and the space of observables is the
algebra
$$
 \mathcal{A} = \overline{P \PsDO^0_{cl}(M;\Lambda^1 M) P},
$$
where $\PsDO^0_{cl}(M;\Lambda^1 M)$ is the algebra of pseudodifferential operators
and $P$ is the orthogonal projection onto the space of co-closed $1$-forms.
The relativistic time evolution is given by the one parameter group generated by
the Laplace-Beltrami operator $\Delta_1$ acting on one-forms.

As we can see the time evolution is in these examples given by a Laplace-type operator acting
on the sections of a vector bundle. The algebra of observables is either the full
algebra of pseudodifferential operators or an appropriate subalgebra that is invariant
under the time-evolution.
The symbol map $\sigma$ is now a map
$$
 \PsDO^0_{cl}(M;E) \to C^\infty(S^* M;\pi^* \mathrm{End}(E))
$$
from the pseudodifferential operators of order $0$ to the smooth
functions on $S^* M$ with values in the bundle $\mathrm{End}(E)$.
Here $\pi^*E$ denotes the pull back of the bundle $E$
on $M$ under the projection $\pi: S^* M \to M$. As in the scalar case the symbol
map has a continuous extension $\hat \sigma$ to a map
$$
 \mathcal{A} \to C(S^* M;\pi^* \mathrm{End}(E))
$$
where $\mathcal{A}=\overline{\PsDO^0_{cl}(M;E)}$ is the norm closure
of the space of pseudodifferential operators acting on $L^2(M;E)$.

In contrast to the Egorov theorem for scalar pseudodifferential operators
the Egorov theorem for matrix valued pseudodifferential operators involves terms
of the time-evolution that are of lower order (\cite{MR661876,EW96}). Let us give some invariant meaning to this.
If $\Delta_E$ is a self-adjoint  Laplace type operator acting on the sections of some hermitian
vector bundle $E$, then there exists a unique connection $\nabla_E$ on $E$ and a unique
potential $V \in C^\infty(M,\mathrm{End}(E))$ such that
$$
 \Delta_E= \nabla_E^* \nabla_E + V.
$$
The locally defined connection-$1$-form can be interpreted as the sub-principal
symbol of $\Delta_E$ (see \cite{JS06}). The connection $\nabla_E$ of course defines a connection 
$\nabla_{\mathrm{End}(E)}$ on $\mathrm{End}(E)$. 
This connection can be used to extend the geodesic
flow on $S^* M$ to a flow on $\pi^*(\mathrm{End}(E))$ by parallel translation. We will denote the
induced action on the sections of $\pi^*(\mathrm{End}(E))$ by $\beta_t$. It is easy to check that
$\beta_t$ is a one-parameter group of $*$-automorphisms on
$$
 C(S^*M; \pi^*(\mathrm{End}(E))).
$$ 
The analog of Egorov's theorem is now as follows. 
Let $\mathcal{A}$ be the algebra $\overline{\PsDO^0_{cl}(M;E)}$.
Then, if $U(t) = e^{-\I t \sqrt{\Delta_E}}$ and $A \in \mathcal{A}$ we have
\begin{gather*}
 A_t=U(-t) A U(t) \in \mathcal{A}, \\
 \hat \sigma_{A_t} = \beta_t(\sigma_{A}).
\end{gather*}
In other words $\beta_t$ is the high energy limit of the quantum time evolution.
A proof can be found in \cite{BO06} and in \cite{JS06}.

\section{Geometric operators and the frame flow}

Most geometric operators like the Dirac operator and the Laplace-Beltrami operator
are acting on sections of vector bundles that are constructed in a geometric way from the
manifold. We assume here that $M$ is oriented and that $FM$ is the bundle of oriented
orthonormal frames. We will show that in many cases the bundle $\mathrm{End}(E)$ as a hermitian
vector bundle with connection is isomorphic to an induced bundle of the frame bundle
by some representation $\rho: SO(n) \to \mathrm{Aut}(\mathrm{End}(\mathbb{C}^m))$
of $SO(n)$ by $*$-automorphisms of $\mathrm{End}(\mathbb{C}^m)$
$$
 \mathrm{End}(E) = FM \times_{\rho} \mathrm{End}(\mathbb{C}^m)
$$
with connection induced by the Levi-Civita connection on $FM$.

\begin{exa}[Dirac operators] Suppose that $D$ is the Dirac operator associated with a spin-
structure or spin$^{c}$-structure acting on the sections of the associated spinor bundle $S$.
Then the action of the complex Clifford algebra bundle $\mathrm{Cl}(TM)$ on $S$
 is irreducible and therefore $\mathrm{End}(S)$ is a quotient of the bundle $\mathrm{Cl}(TM)$.
 The connection on $S$ is compatible with the Clifford action and therefore, the induced connection
 on $\mathrm{End}(S)$ is compatible with the Clifford connection on $\mathrm{Cl}(TM)$.
 But the Clifford algebra bundle is as a hermitian vector bundle with connection obtained
 as an associated bundle
 $$
  \mathrm{Cl}(TM) = FM \times_{\rho} \mathrm{Cl}(\rbb^n),
 $$
 where $\rho$ is the canonical representation of $SO(n)$ on $\mathrm{Cl}(\rbb^n)$.
 Note that the spinor bundle itself is not an associated bundle of $FM$, but $\mathrm{End}(S)$
 nevertheless is. By the Bochner-Lichnerowicz-Weitzenb\"ock-Schr\"odinger formula
 $$
  D^2 = \nabla^* \nabla + V,
 $$
 where $V$ is some potential (for example $\frac{1}{4} R$ in the case of a spin structure).
 Thus, $D^2$ is a Laplace-type operator and the corresponding connection on $S$ and on $\mathrm{End}(S)$
 is the Levi-Civita connection.
 \end{exa}
 
 \begin{exa}
  The bundle $\Lambda^p M$ is an associated bundle of the frame bundle
  $$
   \Lambda^p M = FM \times_{\sigma} \Lambda^p \cbb^n
  $$
  where $\sigma$ is the canonical representation of $SO(n)$ on
  $\Lambda^p \cbb^n$. This of course induces a connection
  on $\Lambda^p M$ which is the Levi-Civita connection on forms.
  The Hodge-Laplace operator on $p$-forms $\Delta_p$ is then defined by
  $$
   \Delta_p = d \delta + \delta d,
  $$
  where $d: C^\infty(\Lambda^p M) \to C^\infty(\Lambda^{p+1} M)$ is the exterior
  differential and $\delta: C^\infty(\Lambda^{p+1} M) \to C^\infty(\Lambda^{p} M)$
  its formal adjoint. Again, by the Bochner--Weitzenb\"ock formula the
  $$
   \Delta_p = \nabla^* \nabla +V,
  $$
  where $V$ involves curvature terms and $\nabla$ is the Levi-Civita connection.
  Note that
  $$
   \mathrm{End}(\Lambda^p M) = FM \times_{\rho} \mathrm{End}(\Lambda^p \cbb^n),
  $$
  where $\rho=\sigma \otimes \sigma^*$.
 \end{exa}
 
These examples show that in geometric situations the bundle $\mathrm{End}(E)$ is
an associated bundle of $FM$. Consequently, the bundle $\pi^* \mathrm{End}(E) \to S^* M$
is an associated bundle of $FM \to S^* M$, where the map $\pi: FM \to S^* M$ is defined by projecting
onto the first vector in the frame and identifying vectors and covectors using the metric.
So, if we view $FM$ as an $SO(n-1)$-principal bundle over $S^*M$ then we can think of
$\pi^* \mathrm{End}(E)$ as the associated bundle
$$
 \pi^* \mathrm{End}(E) = FM \times_{\hat \rho} \mathrm{End}(\cbb^m),
$$
where $\hat \rho$ is the restriction of $\rho: SO(n) \to \mathrm{Aut}(\mathrm{End}(\cbb^m))$
to the subgroup $SO(n-1)$ which we think of as the subgroup that fixes the first vector
in the standard representation on $\rbb^n$.
Therefore, sections of $\pi^* \mathrm{End}(E)$ can be identified with functions $f$ on $FM$
with values in $\mathrm{End}(\cbb^m)$ satisfying the transformation property
\begin{gather} \label{transprop}
 f(x \cdot g) = \hat \rho^{-1}(g) f(x),
\end{gather}
for all $g \in SO(n-1)$ and $x \in FM$.
The frame flow on $FM$ is the extension of the geodesic flow on $S^* M$
by parallel translation to the space $FM$. More explicitly, if $(e_1,\ldots,e_n) \in FM$
is an orthonormal frame, then the frame flow $\gamma_t(e_1,\ldots,e_n)=(e_1(t),\ldots,e_n(t))$
can be defined as follows.
The vector $e_1(t)$ is the tangent of the endpoint of the unique geodesic of length $t$
with starting tangent vector $e_1$. In other words $e_1(t)=G_t(e_1)$.
The rest of the frame $(e_2,\ldots,e_n)$ is parallel transported along this geodesic using the Levi-Civita
connection to give the orthonormal basis $(e_2(t),\ldots,e_n(t))$ in the orthogonal complement of
 $e_1(t)$. It is important here that the Levi-Civita connection preserves angles so that
the parallel transport of the frame yields a frame.
The frame flow gives rise to a flow on the space of functions
on $FM$ by pull-back. By construction $\gamma_t$ commutes with the right action of $SO(n-1)$ on $FM$
and therefore the space of functions satisfying the transformation property (\ref{transprop})
is left invariant. Since all constructions are compatible it turns out that the flow
$\beta_t$ originally constructed from the connection on $E$ coincides with this flow. That is
$$
 (\beta_t f)(x) = f(\gamma_{-t} x),
$$
where $\gamma_t$ is the frame flow.

\section{The high energy limit for geometric operators}

Let $\phi_j$ be a complete  orthonormal sequence of eigensections of a positive Laplace-type operator $\Delta_E$ such that
\begin{gather*}
 \Delta_E \phi_j = \lambda_j^2 \phi_j,\\
 0=\lambda_1 \leq \lambda_2 \leq \lambda_3 \leq \ldots. 
\end{gather*}
Then, the sequence of normal states $\omega_N$ defined by
$$
 \omega_N(A) = \frac{1}{N} \sum_{j=1}^N \langle \phi_j , A \phi_j \rangle
$$
has a high energy limit $\omega_\infty$. As in the scalar case we have
$$
 \omega_\infty(A)=\lim_{t \to 0^+} \frac{\tr (A e^{-t \Delta_E})}{\tr(e^{-t \Delta_E})}
$$
and
$$
  \lim_{t \to 0^+} \frac{\tr (A e^{-t \Delta_E})}{\tr(e^{-t \Delta_E)})} = \frac{1}{\mathrm{rk} E} \int_{S^* } \tr (\hat \sigma_A(\xi)) d\mu_L(\xi),
$$
where $\mu_L$ is the normalized Liouville measure on $S^* M$ and $\tr$ is the trace.
Indeed this state is obviously invariant under the classical time evolution $\beta_t$.

In order to understand which quantum limit can be obtained from subsequences of non-zero counting density
one needs to decompose the state $\omega_\infty$ into ergodic states with respect to the action $\beta_t$.

Indeed, one can prove the following theorem (\cite{MR1384146,JS06,JSZ08})
\begin{theorem}
 Suppose that $p_1,\ldots,p_r$ are projections in $\mathcal{A}$ which commute with $\Delta_E$.
 Suppose furthermore that $\sum_{i=1}^r p_i=\mathrm{id}$ and that the decomposition
 $$
  \omega_\infty(\cdot) = \sum_{i=1}^r \omega_\infty(p_i \cdot)
 $$
 is a decomposition into ergodic states
 $$
  \omega_i(A) = \frac{1}{\omega_\infty(p_i)} \omega_\infty(p_i A).
 $$
 Then Shnirelman's theorem holds in the subspaces onto which $p_i$ project.
 More precisely, if $\phi_j$ is an orthonormal sequence of eigensections of $\Delta_E$ such that
\begin{gather*}
 \Delta_E \phi_j = \lambda_j^2 \phi_j,\\
 \lambda_1 \leq \lambda_2 \leq \lambda_3 \leq \ldots,\\
 p_i \phi_j = \phi_j,
\end{gather*}
and such that $\phi_j$ span the range of $p_i$. Then there is a subsequence 
of eigensections $\phi_{k(j)}$ of counting density one such that
$$
 \lim_{j \to \infty} \langle \phi_{k(j)}, A \phi_{k(j)} \rangle = \omega_i (A).
$$
\end{theorem}

\section{Ergodic decomposition of the tracial state}

Let $\tau : SO(n-1) \to \mathrm{Aut}(\mathrm{End}(\cbb^k))$ be a  representation of $SO(n-1)$
by $*$-automorphisms. Then there is a unique projective unitary representation
$$
 \rho: SO(n-1) \to PU(k)
$$
such that
$$
 \tau(g) (x) = \rho (g)^{-1} x \rho(g).
$$
As before the frame flow induces a flow $\beta_t$ on the space
of sections of 
$$
 F=FM \times_{\tau} \mathrm{End}(\cbb^k)
 $$
 by $*$-automorphisms.
The following theorem is proved in (\cite{JS06,JSZ08}).
\begin{theorem}
 Suppose that the frame flow on $FM$ is ergodic and that $\rho$ is irreducible.
 Then the tracial state
 $$
  \omega(A) = \frac{1}{k} \int_{S^* M } \tr (\hat \sigma_A(\xi)) d\mu_L(\xi),
 $$
 on $C(S^* M;F)$ is ergodic.
\end{theorem}

This gives us a strategy to decompose the tracial state on $C(S^*M, \mathrm{End}(E))$
into ergodic states assuming that $\mathrm{End}(E)$ is an associated bundle of the frame bundle.
Namely, suppose that $\rho: SO(n-1) \to PU(k)$ is a projective unitary representation which gives rise
to a representation $\rho=\rho \otimes \rho^*$ on $\mathrm{End}(\cbb^k)$ such that
$$
  \mathrm{End}(E) = FM \times_{\tau} \mathrm{End}(\cbb^k).
$$
Then decompose $\cbb^k$ into invariant subspaces for $\rho$
$$
 \cbb^k = V_1 \oplus \ldots \oplus V_r.
$$
The orthogonal projection $p_i$ onto $V_i$ is a matrix in $\mathrm{End}(\mathbb{C}^k)$
and we have $\sum_{i=1}^r p_i = \mathrm{id}$. This matrix is invariant under the action
of $\rho$ and therefore, the constant function $p_i \in C^\infty(FM)$
can be understood as a section in
$\mathrm{End}(E)$
as it satisfies the transformation rule.
By the above theorem the state
$$
 \omega_{i}(\cdot) = \frac{1}{\omega(p_i)} \omega(p_i A)
$$
is ergodic and we have constructed an ergodic decomposition of the tracial state.
If there are pseudodifferential operators $P_i$ in $\PsDO^0_{cl}(M,E)$ that are mutually
commuting and commute with $\Delta_E$ such that $\sigma_{P_i} = p_i$
we can then apply the theorem above and conclude that quantum ergodicity holds in subspaces
onto which $P_i$ projects.

\section{Examples}

\subsection{Dirac operators} As before assume that $S$ is the spinor bundle
of some Spin structure or Spin$^c$ structure. Let $D$ be the Dirac operator
and $\mathrm{sign}(D)$ defined by spectral calculus. Then $\mathrm{sign}(D)$ is a 
zero order pseudodifferential operator and its principal symbol is given by
$$
 \sigma_{\mathrm{sign}(D)(\xi)} = \gamma_{\xi},
$$
where $\gamma_\xi$ denotes Clifford multiplication with $\xi$.
Then the projections
$$
 P_{\pm} = \frac{1}{2} (1 \pm \mathrm{sign}(D))
$$
are pseudodifferential operators that commute with $\vert D \vert$.
Their symbols are elements in $C^{\infty}(S^* X, \pi^* \mathrm{End}(S))$
that are invariant under the flow $\beta_t$. If we identify sections of this bundle
with functions on $FM$ with values in the $\mathrm{End}(\cbb^{2^{[n/2]}})$
then this function corresponds to Clifford multiplication with first vector
in the $\mathrm{Cl}(\rbb^n)$-module $\cbb^{2^{[n/2]}}$. This
projects onto an irreducible subspace of the projective representation
$$
 SO(n-1) \to PU(2^{[n/2]}).
$$
If the frame flow is ergodic quantum ergodicity holds in the positive spectral subspace
and negative spectral subspace respectively.

\subsection{Laplace-Beltrami operator on $p$-forms}
 Since the bundles $\Lambda^p T^* M \to M$ are  associated bundles
 of the representation $\Lambda^p \rho$, if $\rho: SO(n) \to \cbb^n$
 the bundle $\pi^* \Lambda^p M \to S^* M$ is associated with 
 the restriction of this representation to the subgroup $SO(n-1)$.
  Note that whereas the representations $\Lambda^p \rho$ are irreducible
  for $p \not= \frac{n}{2}$ the restriction to the group $SO(n-1)$ is not irreducible
  unless $p=0$ or $p=n$. The reason is that since $\cbb^{n}= \cbb \oplus \cbb^{n-1}$
  is a decomposition into invariant subspaces of the $SO(n-1)$-action also
  \begin{gather*}
   \Lambda^p \cbb^n =  \Lambda^p \cbb^{n-1} \oplus \Lambda^{p-1} \cbb^{n-1}
  \end{gather*}
  is a decomposition into invariant subspaces. This decomposition is 
  into irreducible subspaces unless $p = \frac{n-1}{2}$ in which case the first summand is not
  irreducible or $p = \frac{n+1}{2}$ in which case the second summand is not irreducible.
  
  As in the case the Dirac operator one can find pseudodifferential operators of order zero
  that have symbols that project onto these irreducible subspaces.
  Namely, define $\Delta_p^{-1}$ as the inverse of $\Delta_p$ on $(\ker \Delta_p)^{\perp}$
  and to be zero otherwise. Then
  \begin{gather*}
   P= \Delta_p^{-1} \delta d,\\
   Q= \Delta_p^{-1} d \delta .
  \end{gather*}
  are projections  that commute with $\Delta_p$. Their principal symbols are 
 invariant elements in $C^\infty(S^* M, \pi^* \mathrm{End}(\Lambda^p))$
 that give rise to a decomposition of the tracial state as
 $$
  \sigma_P + \sigma_Q = 1.
 $$
 Suppose now that the frame flow on $FM$ is ergodic.
 If $p \not= \frac{n-1}{2}$ then the state
 $$
  \sigma_P(\cdot) = \frac{1}{\omega(P)} \omega(P \cdot)
 $$
 is ergodic as $\sigma_P$ corresponds to the projection onto
 the first summand in.
 If $p \not= \frac{n-1}{2}$ then the state
 $$
  \sigma_Q(\cdot) = \frac{1}{\omega(Q)} \omega(Q \cdot)
 $$
 is ergodic for the analogous reason.
 Thus, quantum ergodicity holds in the subspaces onto which
 $P$ and $Q$ project onto, which are the subspaces of co-exact
 and exact $p$-forms.
 
 In the case $p = \frac{n-1}{2}$ the state $\omega_P$ is not ergodic.
 There is a further pseudodifferential operator commuting with $\Delta_p$
 and with $P$, namely the operator
 $$
  R= \Delta_p^{-1} * d,
 $$
 where $*$ is the Hodge star operator. Note that $R^2=1$ on $\mathrm{rg}(P)$ and the 
 decomposition of the state $\omega_P$ into $+1$ and $-1$ eigenspaces is ergodic in
 case the frame flow is ergodic.
 The eigenvalues of $R$ correspond to polarized forms. For example the case
 of $n=1$, $p=1$ and co-closed $1$-forms corresponds to electrodynamics in dimension
 $3$. It is well known that electromagnetic waves can be decomposed into circular polarized
 waves and  that this decomposition is invariant under the time-evolution determined by
 the Maxwell equation. Sequences of differently polarized eigensections give rise to different quantum limits
 as the observable $R$ that measures the polarization gives rise to an observable in the high energy limit, namely $\sigma_R$,
 that distinguishes them. 
 
 \section{Operators on K\"ahler manifolds}
 The two examples of the previous section may also be discussed in the category of K\"ahler manifolds or other special geometries.
 A K\"ahler manifold can be thought of as a Riemannian manifold of dimension $2m$ such that
 the frame bundle can be reduced to a $U(m)$-principal bundle in such a way that the parallel transport
 preserves the $U(m)$-structure. This is equivalent to the existence of a covariantly constant complex structure.
 On K\"ahler manifolds the frame flow is not ergodic as the complex structure is preserved and gives
 rise to conservation laws. It is much more natural however to consider the $U(m)$-bundle $UM$ of unitary
 frames instead and the restriction of the frame flow to it. Again, a lot of geometric constructions
 in the category of K\"ahler manifolds can be understood as associated bundle constructions
 starting from the unitary frame bundle. The bundles of $(p,q)$-forms are associated bundles of $U M$
 and natural geometric operators to consider are the Dolbeault Dirac operator and the Dolbeault Laplace operator.
 Under the assumption that the unitary frame flow is ergodic the ergodic decomposition of the tracial state
 can be found explicitly. It is closely related to the action of a certain Lie-superalgebra on the space of
 exterior differential forms of a K\"ahler manifold. This action can be seen as the quantum counterpart
 of the classical symmetry that prevents the frame flow on $FM$ to be ergodic. A detailed discussion of this
 and its implications can be found in \cite{JSZ08}.

\section{Conclusions and further remarks}
 It was already found by \cite{MR1644088, MR1694732, MR1838911,MR2073612, BoG04.2} 
 for the Dirac operator in $\mathbb{R}^n$ that the semi-classical limit
 for this operator can be described by a suitable extension of the classical flow and that ergodicity of
 that flow implies quantum ergodicity. The geometric framework discussed in the present article was introduced in
 \cite{JS06} and further developed in \cite{JSZ08}. It deals with high energy limits rather than with the semi-classical
 limit so that only the non-trivial geometry of space contributes to the classical dynamics.
 A representation theoretic framework generalizing the representation theoretic lift on locally symmetric spaces
 to induced bundles over locally symmetric spaces can be found in \cite{BO06}. This article also contains some discussion
 of quantum ergodicity questions for vector bundles.

 As advertised \cite{MR1384146} the language of states and ergodicity of states over $C^*$-dynamical systems 
 is the appropriate one to describe the high energy limit or quantum ergodicity of quantum systems. Its
 application to quantum systems with spin or more mathematically to geometric operators acting on vector bundles
 naturally leads to associated bundle constructions over the frame bundle. The underlying dynamics being the frame flow.
 Quantum ergodicity questions translate into questions about the frame flow. In particular ergodicity of the frame
 flow has strong implications for quantum systems with spin. It implies quantum ergodicity on certain natural subspaces
 that can be found more or less constructively from our method. Finally we would like to mention that
 the frame flow was already considered by Arnold in \cite{MR0158330}.  In
negative curvature, it was studied by Brin, together with Gromov,
Karcher and Pesin, in a series of papers \cite{MR0343316, MR0370660,MR0394764, MR0582702, MR0670078, MR0756723}
independently of any connection to spectral theory for operators on vector bundles or quantum ergodicity questions.
Many examples of manifolds with ergodic frame flow are known (for example manifolds of constant negative curvature to name only the simplest ones) and much progress has been made towards its understanding. We would like to refer the reader to the above mentioned literature on
frame flows for further details.

\vspace{1cm}

\noindent{\bf Acknowledgements.} The author would like to thank the CRM Montreal and the Analysis lab in Montreal
kind hospitality during his stay in summer 2008

\end{document}